\documentclass[12pt]{amsart}
\newcommand{\RP}{{\mathbb {RP}}}
\newcommand{\CP}{{\mathbb {CP}}}
\newcommand{\Z}{{\mathbb Z}}
\newcommand{\R}{{\mathbb R}}
\newcommand{\C}{{\mathbb C}}
\newcommand{\Q}{{\mathbb Q}}
\newcommand{\sm}{\setminus}

\newtheorem{theorem}{Theorem}
\newtheorem{example}{Example}
\newtheorem{proposition}{Proposition}
\newtheorem{corollary}{Corollary}
\newtheorem{remark}{Remark}
\newtheorem{definition}{Definition}
\newtheorem{lemma}{Lemma}

\begin{document}

\author{V.A.~Vassiliev}
\address{Steklov Mathematical Institute of Russian Academy of Sciences;
National Research Institute --- Higher School of Economics}
\email{vva@mi.ras.ru}


\thanks{
Supported by RSF grant 14-50-00005}

\title{Homology groups of spaces of non-resultant quadratic polynomial systems in
$\R^3$} 

\begin{abstract}
We calculate rational homology groups of spaces of non-resultant (i.e. having
no non-trivial common zeros) systems of homogeneous quadratic polynomials in
$\R^3$.
\end{abstract}

\maketitle

\section{Introduction}

The {\em resultant subvariety} in the space of systems of homogeneous
polynomials in $\R^n$ consists of systems having common zeros besides the
origin in $\R^n$. The varieties of this kind were considered, in particular, by
S.~Smale and his collaborators in connection with the $P \stackrel{?}{=} NP$
problem over fields of zero characteristic, see \cite{Smale}, \cite{BSS},
\cite{SS}. Also, the topology of complements of resultant varieties in certain
spaces of systems of non-homogeneous polynomials in one complex variable was
studied in \cite{CCMM}, \cite{Kozl}, \cite{V}. The homology groups of spaces of
non-resultant homogeneous systems in $\R^2$ were calculated in \cite{homres}.
Below we study the next complicated problem --- of homogeneous quadratic
polynomial systems in $\R^3$.

Let $W \simeq \R^6$ be the space of all real quadratic forms in $\R^3$. For any
natural $k$, consider the space $W^k \simeq \R^{6k}$ of all quadratic maps
$\R^3 \to \R^k$. Denote by $\Sigma$ the {\em resultant subvariety} in $W^k$,
i.e. the space of all systems of $k$ quadratic forms having common zeros in
$\R^3 \setminus 0$.

\begin{theorem}
\label{mthm12} For any even $ k \ge 4$, the Poincar\'e polynomial of the group
$H^*(W^k \setminus \Sigma, \Q)$ is equal to
\begin{equation} \label{stabeven}(1+t^{k-1})(1+t^{3k-9}+t^{5k-14}). \end{equation}

For any odd $ k \ge 3$, the Poincar\'e polynomial of $H^*(W^k \setminus \Sigma,
\Q)$ is equal to
\begin{equation} \label{stabodd}1+t^{k-1} + t^{3k-7}(1+t^{k-5})(1+t^{2k-3}).
\end{equation}

The group $H^i(W^2 \setminus \Sigma, \Q)$ is equal to $\Q$ if $i=0$ or $i=1$,
and is trivial otherwise.
\end{theorem}

\begin{example}
\label{exam3} \rm For $k=3$, the polynomial (\ref{stabodd}) is equal to
$(2+t^3)(1+t^2)$. In particular, the space $W^3 \setminus \Sigma$ consists of
two connected components. Indeed, any non-resultant triple of quadratic forms
in $\R^3$ defines a map $\RP^2 \to \R^3 \setminus 0 \sim S^2.$ The $ \mbox{mod
}2$ degree of this map is invariant along any connected component of $W^3
\setminus \Sigma$. This degree is equal to 0 for any system containing
non-negative quadratic forms, and is non-trivial for a triple of quadratic
forms whose zero sets in $\RP^2$ meet like
 \unitlength=0.3mm \special{em:linewidth 0.4pt} \linethickness{0.4pt}
\begin{picture}(30,27)
\put(10,10){\circle{20}} \put(20,10){\circle{20}} \put(15,17){\circle{20}}
\end{picture}.
\end{example}

\section{First reductions}

\subsection{Alexander duality}

Denote by $\overline{H}_*$ the {\em Borel--Moore homology group}, i.e. the
homology group of the complex of locally finite singular chains of a
topological space. In the case of homology groups with coefficients in a group
(i.e. in a constant local system) $A$, an equivalent definition is as follows:
$\overline{H}_*(X,A)$ is the homology group of the one-point compactification
of $X$ reduced modulo the added point. However, in the case of non-constant
coefficient systems the latter definition is generally not applicable. $\tilde
H_*$ and $\tilde H^*$ is the notation for the homology and cohomology groups
reduced modulo a point.
\medskip

Following \cite{A70}, we use the Alexander isomorphism
\begin{equation}
\label{alex} \tilde H^i(W^k \sm \Sigma) \simeq \overline{H}_{6k-i-1}(\Sigma).
\end{equation}

\subsection{A toy calculation: homology of Stiefel manifolds}
\label{toy}

To calculate the right-hand group in (\ref{alex}) we consider a {\it simplicial
resolution} of the space $\Sigma$. Let us demonstrate this method in a simpler
situation, replacing quadratic maps by linear ones. The space of non-resultant
collections of $k$ linear functions $\R^3 \to \R^1$, $k \ge 3$, is homotopy
equivalent to the Stiefel manifold $V_3(\R^k)$. Its cohomology group is
well-known; let us show how to calculate it by our techniques.

\begin{proposition} \label{toy12}
For any $k \ge 3$, the Poincar\'e polynomial of the group $H^i(V_3(\R^k), \Q)$
is equal to $(1+t^{k-1})(1+t^{2k-5})$ if $k$ is even, and to
$(1+t^{k-3})(1+t^{2k-3})$ if $k$ is odd.
\end{proposition}

{\it Proof.} Consider the entire space $\R^{3k}$ of collections of $k$ linear
functions, and the discriminant variety $\Sigma$ in it consisting of forms
vanishing at some non-trivial subspace. To construct the simplicial
resolution $\sigma$ of this variety, we list all subsets in $\RP^2$
which can be distinguished by systems of linear equations.

It are 1) particular points (the space of them is $\RP^2$ itself), 2) the lines
(the space of them is the dual projective plane $\widehat{\RP^2}$), and 3)
entire $\RP^2$: this is a single object named $[\RP^2]$. Let us denote these
spaces of subsets by $J_1$, $J_2$ and $J_3$ respectively.

We consider them as topological spaces with their standard topology. The
disjoint union $J_1 \sqcup J_2 \sqcup J_3$ is a partially ordered set with the
order defined by inclusion of corresponding subsets. Consider the topological
order complex $\Psi$ of this poset. Namely, take the join $J_1
* J_2 * J_3$ of these three spaces (i.e., the naturally topologized union of
triangles and segments, whose vertices correspond to the points of different
spaces) and define the space $\Psi$ as the union of only those of these
simplices, all whose vertices are incident to one another. For any point $X$ of
$J_1, J_2$ or $J_3$, define its subordinate order subcomplex $\Psi(X)$ as the
union of all such simplices, all whose vertices correspond to sets contained in
the set $\{X\}$. Also, denote by $\partial \Psi(X)$ the {\em link} of
$\Psi(X)$, i.e. the union of all simplices participating in $\Psi(X)$ but not
containing its vertex $X$. Let $\breve \Psi(X)$ be the difference $\Psi(X)
\setminus
\partial \Psi(X)$. Denote by $\Psi_p$, $p=1, 2, 3$ the union of all complexes
$\Psi(X)$ over all $X \in J_j, $ $j \le p$. Obviously, $\Psi_3 = \Psi([\RP^2])
= \Psi$.

The {\em simplicial resolution} $\sigma$ of $\Sigma$ is defined as a subspace
in the direct product $\Psi \times \R^{3k}$. Namely, for any $p=1,2,$ or $3$,
and any $X \in J_p$ we take the product $\Psi(X) \times L(X)$, where $L(X)
\subset \R^{3k}$ is the space of all systems of $k$ linear functions vanishing
on $X$. $\sigma$ is defined as the union of all such products over all $p$ and
all $X \in J_p$. The subspace $\sigma_p \subset \sigma$ is the union of such
products only over all $X \in J_j,$ $j \le p$. It is easy to see that the
difference $\sigma_p \setminus \sigma_{p-1}$ is the space of a fiber bundle
over $J_p$; its fiber over $X \in J_p$ is equal to $\breve \Psi(X) \times
L(X)$. In particular, $\sigma_p \setminus \sigma_{p-1}$ is the space of a
vector bundle over $\Psi_p \setminus \Psi_{p-1}$ with the fibers equal to
$L(X)$, $X \in J_p$.

The projection $\sigma \to \Sigma$ is proper with contractible fibers, hence
defines an isomorphism $\overline{H}_*(\sigma, \Q) \to
\overline{H}_*(\Sigma,\Q)$. Now let us calculate the spectral sequence
converging to the group $\overline{H}_*(\sigma, \Q)$ and  generated by the
filtration $\sigma_1 \subset \sigma_2 \subset \sigma_3 \equiv \sigma$. By
definition, its initial term $E^1_{p,q}$ is equal to
$\overline{H}_{p+q}(\sigma_p \setminus \sigma_{p-1}, \Q)$.

The space $\sigma_1$ is fibered over $J_1 = \RP^2$. Its fiber over any point $x
\in \RP^2$ is isomorphic to $\R^{2k}$ and consists of all linear systems
vanishing on the corresponding line in $\R^3$. This fiber bundle is orientable
if and only if $k$ is even, so by Thom isomorphism we have
$\overline{H}_i(\sigma_1, \Q) \simeq H_{i-2k}(\RP^2,\Q)$ if $k$ is even, and
$\overline{H}_i(\sigma_1, \Q) \simeq H_{i-2k}(\RP^2,\Q \otimes {\mathbb O}r)$
if $k$ is odd; here ${\mathbb O}r$ is the orientation bundle of $\RP^2$.

The space $\sigma_2 \setminus \sigma_1$ is fibered over $J_2=\widehat{\RP^2}$.
Its fiber over any point $l \in \widehat{\RP^2}$ is the product of the
$k$-dimensional subspace $ L(l) \subset \R^{3k}$ and the open two-dimensional
disc $\breve \Psi(l)$, which is the union of all half-open intervals connecting
the point $l \in J_2$ with all points of $J_1$ corresponding to the points of
the line $l$. The fiber bundle over $J_2$ with fibers $L(l)$ is orientable if
and only if $k$ is even; the fiber bundle with fibers $\breve \Psi(l)$ always
is non-orientable. Therefore $\overline{H}_i(\sigma_2 \setminus \sigma_1, \Q)
\simeq H_{i-k-2}(\RP^2, \Q)$ if $k$ is odd, and $\overline{H}_i(\sigma_2
\setminus \sigma_1, \Q) \simeq H_{i-k-2}(\RP^2, \Q \otimes {\mathbb O}r)$ if
$k$ is even. Finally, $\sigma_3 \setminus \sigma_2$ is in all cases the entire
order complex $\Psi$ of non-empty linear subspaces of $\R^3$, from which its
order subcomplex $\partial \Psi$ of only proper subspaces is removed. $\Psi$ is
a cone over this subcomplex, which is known (see e.g. \cite{BS}, \cite{Vfil})
to be homeomorphic to $S^4$. So, $\overline{H}_i(\sigma_3 \setminus \sigma_2,
\Q)$ is equal to $\Q$ if $i=5$, and is trivial otherwise. Summarizing, we get
that our spectral sequence has only the following non-trivial terms $E^1_{p,q}$
(all of which are equal to $\Q$). For even $k$ this are terms with $(p,q)$
equal to $(1, 2k-1), (2, k+2),$ and $(3,2)$; for odd $k$ the terms with $(p,q)$
equal to $(1,2k+1), (2,k),$ and $(3,2)$ only. All further differentials of the
spectral sequence are obviously trivial. We obtain the group
$\overline{H}_*(\sigma, \Q)$, and Proposition \ref{toy12} follows by Alexander
duality similar to (\ref{alex}) but with $6k$ replaced by $3k$. \hfill $\Box$

\subsection{Singular sets for the resultant in the spaces of quadratic maps}
\label{list}

From now on, we come back to the space $W^k$ of systems of $k$ quadratic
functions $\R^3 \to \R$, consider the resultant subvariety $\Sigma$ in it, and
construct the simplicial resolution $\sigma$ of it. In this construction, we
follow the A.~Gorinov's (see \cite{Gorinov}, \cite{Tommasi2}) modification of
the construction of resolutions of discriminant sets proposed in
\cite{novikov}.

We fix the following list of families of compact subsets in $\RP^2$, basic for
the construction of the resolution, cf. \cite{novikov}, \cite{Gorinov}.

{\bf $J_1 \sim \RP^2$} consists of particular points.

{\bf $J_2$} is the configuration space $B(\RP^2,2)$ of unordered pairs of
distinct points in $\RP^2$.

{\bf $J_3$} is the similar space $B(\RP^2,3)$ of subsets of cardinality 3 in
$\RP^2$.

{\bf $J_4$} is the space of quadruples of distinct points, all of which belong
to one and the same line in $\RP^2$.

{\bf $J_5$} is the space $\widehat{\RP^2}$ of lines.

{\bf $J_6$} is the space $B^\times(\RP^2,4)$ of quadruples of different points
in $\RP^2$, which do not belong to one and the same line.

{\bf $J_7$} is the space of configurations of the form $\{$a line, a point
outside it in $\RP^2\}$.

{\bf $J_8$} is the space $B(\widehat{\RP^2},2)$ of pairs of different lines in
$\RP^2$.

{\bf $J_9$} is the space of non-singular non-empty conics in $\RP^2$.

{\bf $J_{10}$} consists of unique point and corresponds to entire $\RP^2$.

\medskip
Following \cite{Gorinov}, supply the union of these spaces with the Hausdorff
metric on the space of compact subsets in $\RP^2$: the distance between two
compact subsets $X, Y \subset \RP^2$ is
\begin{equation}
\label{hausd} \max_{x \in X} \rho(x, Y) + \max_{y \in Y} \rho(y, X),
\end{equation}
where $\rho$ is the usual metric in the plane $\RP^2$ considered as a quotient
space of the unit 2-sphere. Denote by ${\mathfrak J}$ the obtained metric space
$J_1 \cup J_2 \cup \dots \cup J_{10}$. It is easy to check that it is a closed
subset of the space of all compact subsets of $\RP^2$ supplied with the metric
(\ref{hausd}). The latter space is compact, hence our space ${\mathfrak J}$ is
compact too.

Again, consider this space ${\mathfrak J}$ as a partially ordered set by
incidence of corresponding subsets. Unlike the situation considered in \S
\ref{toy}, some of its strata $J_p$ are not compact and adjoin one another. By
this reason the construction of the simplicial resolution is slightly more
complicated. In particular, we use the following notion (see \cite{BC},
\cite{VO}).
\medskip

\begin{definition} \rm
Given a compact finite-dimensional topological space $X$, let us embed it
generically into a space $\R^\omega$ of a very large dimension, and define its
{\it $r$-th self-join} $X^{*r}$ as the union of all simplices of dimensions $1,
\dots, r-1$ in $\R^\omega$, all whose vertices belong to the image of $X$. (The
genericity of an embedding means that intersections of different simplices are
their common faces only).
\end{definition}

It is easy to see that the spaces $X^{*r}$ defined by different generic
embeddings of $X$ are canonically homeomorphic to one another.

Define the order complex $\Psi$ of sets $J_1, \dots, J_{10}$ as a subspace of
the 10-th self-join ${\mathfrak J}^{*10}$: it consists of all possible
simplices (of different dimensions) in  ${\mathfrak J}^{*10}$, all whose
vertices are incident to one another. For any $X \in {\mathfrak J}$, define the
space $\Psi(X)$ as the union of simplices in $\Psi$, all whose vertices
correspond to sets containing $\{X\}$. Also we define the link $\partial
\Psi(X)$ as the union of simplices in $\Psi(X)$ not containing the point $X$,
and the space $\breve \Psi(X)$ as the difference $\Psi(X) \setminus \partial
\Psi(X)$. By the exact sequence of the pair we have
$$\overline{H}_i(\breve \Psi(X)) \equiv \tilde H_{i-1}(\partial \Psi(X)).$$

The space $\Psi$ is naturally filtered: for any $p \in [1, 10]$ the space
$\Psi_p \subset \Psi$ is defined as the union of spaces $\Psi(X)$ over all $X
\in J_j$, $j \le p$. By the construction, $\Psi_p \setminus \Psi_{p-1}$ is the
union of spaces $\breve \Psi(X)$ over all $X \in J_p$.

For any point $X \in J_p$ define the space $\sigma(X) \subset \Psi \times W^k$
as the direct product $\Psi(X) \times L(X)$, where $L(X)$ is the space of all
systems vanishing on the set $\{X\}$. The space $\sigma_p$ is defined as the
union of spaces $\sigma(X)$ over all $X \in J_j,$ $j \le p$; $\sigma \equiv
\sigma_{10}$.  So, we have filtrations
\begin{equation}
\Psi_1 \subset \dots \subset  \Psi_{10} \equiv  \Psi, \hspace{2cm} \sigma_1
\subset \dots \subset  \sigma_{10} \equiv \sigma.
 \label{fltr2}
\end{equation}

The obvious projection $ \Psi \times W^k \to W^k$ induces a map $\Pi: \sigma
\to \Sigma.$ For any point $F \in \Sigma$ (i.e. some resultant system of
polynomials) its preimage $\Pi^{-1}(F)$ is the subcomplex $ \Psi(X(F)) \times
F$, where $X(F) \in {\mathfrak J}$ is the whole set of common zeros of all
polynomials from $F$. (Any such set belongs to some $J_p$).

\begin{lemma} \label{lem2}
1. The map $\Pi: \sigma \to \Sigma$ induced by the obvious projection $ \Psi
\times W^k \to W^k$ is proper, and its continuation to a map of one-point
compactifications of these spaces $\sigma$ and $\Sigma$ is a homotopy
equivalence. In particular, it induces an isomorphism of Borel--Moore homology
groups, $\overline{H}_*(\sigma) \simeq \overline{H}_*(\Sigma)$.

2. For any $p=1, \dots, 10$, the projection $\sigma \to \Psi$ maps $\sigma_p$
to $\Psi_p$. Its restriction to $\sigma_p \setminus \sigma_{p-1}$ is a vector
bundle over $\Psi_{p} \setminus \Psi_{p-1}$, whose dimension is equal to that
of the space $L(X)$ for any $X \in J_p$.

3. For any $p$, $\Psi_p \setminus \Psi_{p-1}$ is locally trivially fibered over
$J_p$, with the fiber over $X \in J_p$ equal to $\breve \Psi(X)$.
\end{lemma}

For the proof of item 1, see \cite{Gorinov}, Theorem 2.8. Indeed, our
collection of families $J_1, \dots, J_{10}$ obviously satisfies all conditions
listed in \cite{Gorinov}, p. 399 (see also \cite{Tommasi2}, pp. 13--14). Items
2 and 3 follow immediately from the construction. \quad $\Box$
\medskip

So, we can (and will) calculate the group $\overline{H}_*(\sigma)$ instead of
$\overline{H}_*(\Sigma)$.

Consider the spectral sequence $E_{p,q}^r,$ calculating the group
$\overline{H}_*(\sigma, \Q)$ and generated by the filtration (\ref{fltr2}),
with unique reservation: we do not distinguish $\sigma_4$ from the rest of
$\sigma_5$, so that the revised filtration looks as follows: $\sigma_1 \subset
\sigma_2 \subset \sigma_3 \subset \sigma_5 \subset \dots$. The term $E_{p,q}^1$
is thus canonically isomorphic to the group $\overline{H}_{p+q}(\sigma_p \sm
\sigma_{p-1}, \Q)$ if $p=1, 2$ or $3$, to $\overline{H}_{4+q} (\sigma_5
\setminus \sigma_3, \Q)$ if $p=4$, and $\overline{H}_{p+q}(\sigma_{p+1} \sm
\sigma_p, \Q)$ for $p=5, \dots , 9$.

\unitlength=0.70mm \special{em:linewidth 0.4pt} \linethickness{0.4pt} \noindent
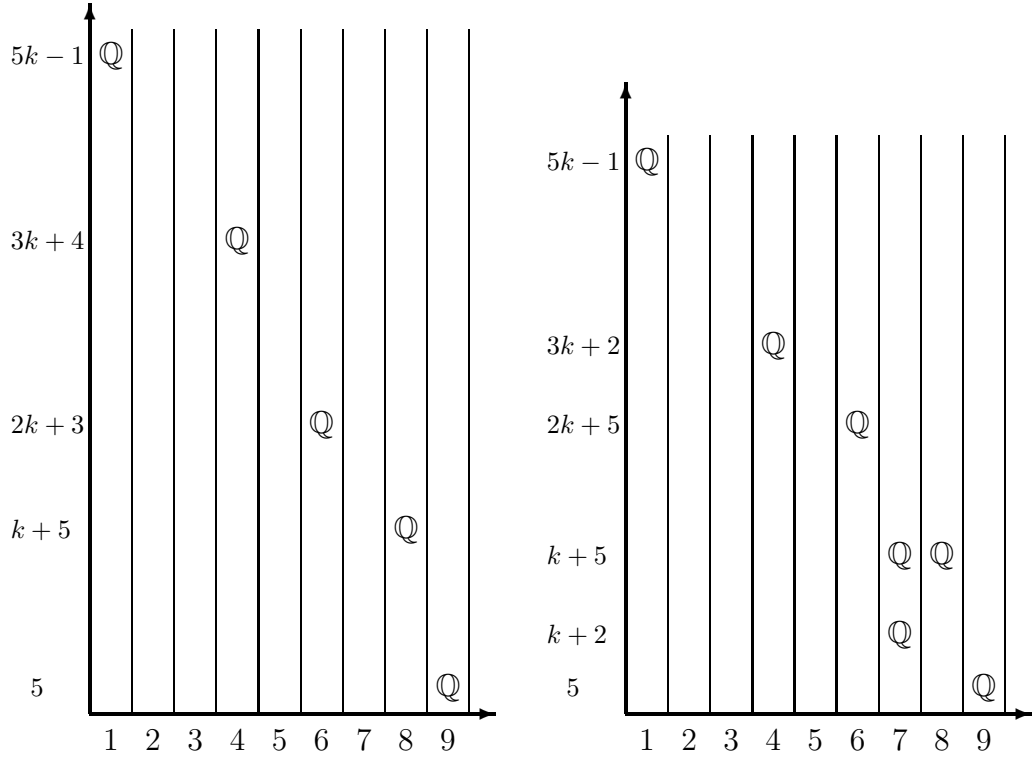
\begin{figure}
\mbox{
\begin{picture}(90.00,145.00)
\thicklines \put(15.00,10.00){\vector(1,0){77.00}}
\put(15.00,10.00){\vector(0,1){135.00}} \thinlines
\put(23.00,10.00){\line(0,1){130.00}} \put(31.00,10.00){\line(0,1){130.00}}
\put(39.00,10.00){\line(0,1){130.00}} \put(47.00,10.00){\line(0,1){130.00}}
\put(55.00,10.00){\line(0,1){130.00}} \put(63.00,10.00){\line(0,1){130.00}}
\put(71.00,10.00){\line(0,1){130.00}} \put(79.00,10.00){\line(0,1){130.00}}
\put(87.00,10.00){\line(0,1){130.00}}
\put(0.00,135.00){\makebox(0,0)[lc]{\footnotesize $5k-1$}}
\put(19.00,135.00){\makebox(0,0)[cc]{$\Q$}}
\put(0.00,100.00){\makebox(0,0)[lc]{\footnotesize $3k+4$}}
\put(43.00,100.00){\makebox(0,0)[cc]{$\Q$}}
\put(0.00,65.00){\makebox(0,0)[lc]{\footnotesize $2k+3$}}
\put(59.00,65.00){\makebox(0,0)[cc]{$\Q$}}
\put(0.00,45.00){\makebox(0,0)[lc]{\footnotesize $k+5$}}
\put(75.00,45.00){\makebox(0,0)[cc]{$\Q$}}
\put(5.00,15.00){\makebox(0,0)[cc]{\footnotesize $5$}}
\put(83.00,15.00){\makebox(0,0)[cc]{$\Q$}}
\put(19.00,5.00){\makebox(0,0)[cc]{1}} \put(27.00,5.00){\makebox(0,0)[cc]{2}}
\put(35.00,5.00){\makebox(0,0)[cc]{3}} \put(43.00,5.00){\makebox(0,0)[cc]{4}}
\put(51.00,5.00){\makebox(0,0)[cc]{5}} \put(59.00,5.00){\makebox(0,0)[cc]{6}}
\put(67.00,5.00){\makebox(0,0)[cc]{7}} \put(75.00,5.00){\makebox(0,0)[cc]{8}}
\put(83.00,5.00){\makebox(0,0)[cc]{9}}
\end{picture}
} \quad \mbox{
\begin{picture}(68.00,130.00)
\thicklines \put(15.00,10.00){\vector(1,0){77.00}}
\put(15.00,10.00){\vector(0,1){120.00}} \thinlines
\put(23.00,10.00){\line(0,1){110.00}} \put(31.00,10.00){\line(0,1){110.00}}
\put(39.00,10.00){\line(0,1){110.00}} \put(47.00,10.00){\line(0,1){110.00}}
\put(55.00,10.00){\line(0,1){110.00}} \put(63.00,10.00){\line(0,1){110.00}}
\put(71.00,10.00){\line(0,1){110.00}} \put(79.00,10.00){\line(0,1){110.00}}
\put(87.00,10.00){\line(0,1){110.00}} \put(19.00,5.00){\makebox(0,0)[cc]{1}}
\put(27.00,5.00){\makebox(0,0)[cc]{2}} \put(35.00,5.00){\makebox(0,0)[cc]{3}}
\put(43.00,5.00){\makebox(0,0)[cc]{4}} \put(51.00,5.00){\makebox(0,0)[cc]{5}}
\put(59.00,5.00){\makebox(0,0)[cc]{6}} \put(67.00,5.00){\makebox(0,0)[cc]{7}}
\put(75.00,5.00){\makebox(0,0)[cc]{8}} \put(83.00,5.00){\makebox(0,0)[cc]{9}}
\put(0.00,115.00){\makebox(0,0)[lc]{\footnotesize $5k-1$}}
\put(19.00,115.00){\makebox(0,0)[cc]{$\Q$}}
\put(0.00,80.00){\makebox(0,0)[lc]{\footnotesize $3k+2$}}
\put(43.00,80.00){\makebox(0,0)[cc]{$\Q$}}
\put(0.00,65.00){\makebox(0,0)[lc]{\footnotesize $2k+5$}}
\put(59.00,65.00){\makebox(0,0)[cc]{$\Q$}}
\put(0.00,40.00){\makebox(0,0)[lc]{\footnotesize $k+5$}}
\put(75.00,40.00){\makebox(0,0)[cc]{$\Q$}}
\put(67.00,40.00){\makebox(0,0)[cc]{$\Q$}}
\put(0.00,25.00){\makebox(0,0)[lc]{\footnotesize $k+2$}}
\put(67.00,25.00){\makebox(0,0)[cc]{$\Q$}}
\put(5.00,15.00){\makebox(0,0)[cc]{\footnotesize $5$}}
\put(83.00,15.00){\makebox(0,0)[cc]{$\Q$}}
\end{picture}
} \caption{$E^1$ for $k$ even (k=6) \hspace{1.2cm} and $k$ odd (k=5)}
\label{mss}
\end{figure}

\begin{theorem}
\label{prop21} For any even $($respectively, odd$)$ $k$, all non-trivial groups
$E^1_{p,q}$ of the initial term of our spectral sequence are as shown in Fig.
\ref{mss} left $($respectively, right$)$. In the case of odd $k$, the
differential $d^1: E^1_{8,k+5} \to E^1_{7,k+5}$ is an isomorphism.
\end{theorem}

Theorem \ref{mthm12} for $k \ge 3$ follows immediately from this one, because
no other differentials $d^r: E^r_{p,q} \to E^r_{p-r, q+r-1}$ of the spectral
sequence can connect non-trivial cells of these tables.

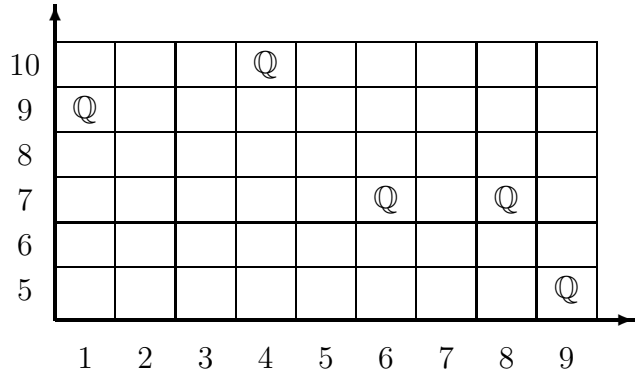
\begin{figure}
\unitlength=1mm
\begin{picture}(87.00,52.00)
\thicklines \put(10.00,10.00){\vector(1,0){77.00}}
\put(10.00,10.00){\vector(0,1){42.00}} \thinlines
\put(18.00,10.00){\line(0,1){37.00}} \put(26.00,10.00){\line(0,1){37.00}}
\put(34.00,10.00){\line(0,1){37.00}} \put(42.00,10.00){\line(0,1){37.00}}
\put(50.00,10.00){\line(0,1){37.00}} \put(58.00,10.00){\line(0,1){37.00}}
\put(66.00,10.00){\line(0,1){37.00}} \put(74.00,10.00){\line(0,1){37.00}}
\put(82.00,10.00){\line(0,1){37.00}} \put(10.00,17.00){\line(1,0){72.00}}
\put(10.00,23.00){\line(1,0){72.00}} \put(10.00,29.00){\line(1,0){72.00}}
\put(10.00,35.00){\line(1,0){72.00}} \put(10.00,41.00){\line(1,0){72.00}}
\put(10.00,47.00){\line(1,0){72.00}} \put(6.00,14.00){\makebox(0,0)[cc]{$5$}}
\put(6.00,20.00){\makebox(0,0)[cc]{$6$}}
\put(6.00,26.00){\makebox(0,0)[cc]{$7$}}
\put(6.00,32.00){\makebox(0,0)[cc]{$8$}}
\put(6.00,38.00){\makebox(0,0)[cc]{$9$}}
\put(6.00,44.00){\makebox(0,0)[cc]{$10$}}
\put(70.00,26.00){\makebox(0,0)[cc]{$\Q$}}
\put(78.00,14.00){\makebox(0,0)[cc]{$\Q$}}
\put(54.00,26.00){\makebox(0,0)[cc]{$\Q$}}
\put(38.00,44.00){\makebox(0,0)[cc]{$\Q$}}
\put(14.00,38.00){\makebox(0,0)[cc]{$\Q$}}
\put(14.00,5.00){\makebox(0,0)[cc]{1}} \put(22.00,5.00){\makebox(0,0)[cc]{2}}
\put(30.00,5.00){\makebox(0,0)[cc]{3}} \put(38.00,5.00){\makebox(0,0)[cc]{4}}
\put(46.00,5.00){\makebox(0,0)[cc]{5}} \put(54.00,5.00){\makebox(0,0)[cc]{6}}
\put(62.00,5.00){\makebox(0,0)[cc]{7}} \put(70.00,5.00){\makebox(0,0)[cc]{8}}
\put(78.00,5.00){\makebox(0,0)[cc]{9}}
\end{picture}
\caption{$E^1$ for $k=2$} \label{sstwo}
\end{figure}

For $k=2$ the left-hand table of Fig. \ref{mss} becomes as is shown in Fig.
\ref{sstwo}. Since the entire space $W^2$ is 12-dimensional, all groups
$E_{p,q}$ with $p+q \ge 12$ should be killed by some differentials. This can
happen in unique way: the differentials $d^4: E^4_{8,7} \to E^4_{4,10}$ and
$d^3: E^3_{9,5} \to E^3_{6,7}$ should be isomorphisms. Thus, the only surviving
term is $E_{1,9} \simeq \Q$. This gives us $\overline{H}_{10}(\Sigma,\Q) \simeq
\Q$ and proves (via the Alexander duality) the assertion of Theorem
\ref{mthm12} for the remaining case $k=2$.
\medskip

The proof of Theorem \ref{prop21} takes the rest of this article.

\begin{remark} \rm
In both spectral sequences shown in Fig. \ref{mss}, the column $p=9$ consists
of the Borel--Moore homology groups of $\sigma_{10} \setminus \sigma_9 \equiv
\Psi_{10} \setminus \Psi_9 \equiv \breve \Psi([\RP^2]).$ Namely, $E^1_{9,q}
\simeq \overline{H}_{9+q}(\breve \Psi([\RP^2]), \Q) \equiv \tilde
H_{8+q}(\partial \Psi([\RP^2]),\Q)$. The latter space $\partial \Psi([\RP^2])$
is the simplicial resolution of the topological order complex of the partially
ordered set of all proper subsets in $\RP^2$, listed in \S \ref{list}. Thus we
prove in particular that this order complex has rational homology group of
$S^{13}$. This order complex is very similar (but does not coincide) with the
order complex of all proper subspaces in $W$ consisting of quadrics vanishing
on some subset in $\RP^2$. The latter complex is studied in more detail in
\cite{orcom}. In particular, it is proved there that this complex also has
rational homology of $S^{13}$.

Notice that the similar order complex of sets defined by systems of linear
functions in $\R^3$ is homeomorphic to $S^4$; we have used this fact in the
proof of Proposition \ref{toy12}.
\end{remark}

\section{Preliminary facts on homology of configuration spaces}

Given a topological space $M$, denote by $B(M,j)$ the {\em configuration
space}, whose points are the subsets of cardinality $j$ in $M$.

\begin{lemma} \label{lem1} For any $j,$ there is a locally trivial
fibre bundle $B(S^1,j) \to S^1$, whose fiber is homeomorphic to $\R^{j-1}$.
This fibre bundle is orientable $($and hence trivial$)$ if $j$ is odd and is
non-orientable if $j$ is even. \end{lemma}

Indeed, consider $S^1$ as the unit circle in $\C^1$, then the projection map of
this fibre bundle can be realised as taking the product of $j$ complex numbers
on this circle. The fibre of this bundle can be identified in the terms of the
universal covering $\R^j \to T^j$ with any connected component of some
hyperplane $\{x_1 + \dots + x_j = \mbox{const}\} \subset \R^j,$ from which all
affine planes given by $x_m = x_n + 2\pi k$, $m \neq n$, $k \in \Z$, are
removed. Any such component is convex and hence diffeomorphic to $\R^{j-1}$.
The assertion on the orientability can be checked immediately. \hfill $\Box$
\medskip

\begin{remark} \label{rem00} \rm The image of the standard section of this
fiber bundle consists of all $j$-ples of points separating $S^1$ into equal
arcs.
\end{remark}

For any manifold $M$, we denote by ${{\mathbb O}r}$ the orientation local
system of groups on $M$: it is locally isomorphic to the constant $\Z$-system,
but the elements of $\pi_1(M)$ violating the orientation of $M$ act on it as
the multiplication by $-1$. For any topological space $M$ denote by $\pm \Z$
the local system on $B(M,j)$, which is locally isomorphic to the constant
$\Z$-bundle, but the elements of $\pi_1(B(M,j))$ defining odd permutations of
$j$ points act on it as multiplication by $-1$. Also, set $\pm \Q \equiv \pm \Z
\otimes \Q$.


\begin{lemma}
\label{lem77} For any even-dimensional manifold $M$, a closed loop in the
manifold $B(M,j)$ violates the orientation of this manifold if and only if the
union of traces of all $j$ points of our configuration during their
corresponding movement defines a disorienting cycle in $H_1(M, \Z_2)$ $($i.e. a
cycle taking non-zero value on the first Stiefel--Whitney class of $M)$.
\end{lemma}

{\it Proof} is trivial, see e.g. \cite{V}. \hfill $\Box$
\medskip

Consider the $j!$-fold covering $\nu: I(\RP^2,j) \to B(\RP^2,j)$, where
$I(\RP^2,j)$ is the space of {\em ordered} sets of $j$ distinct points in
$\RP^2$.

\begin{lemma} \label{lemord}
For any natural $j \ge 2$ the Borel--Moore homology group
$\overline{H}_*(I(\RP^2,j), \Q)$ is trivial in all dimensions.
\end{lemma}

{\it Proof.} $I(\RP^2,j) = (\RP^2)^j \setminus \Delta$, where $\Delta \subset
(\RP^2)^j$ is the set of all sequences of $j$ points in $\RP^2$, at least two
of which coincide. By the exact sequence of the pair $((\RP^2)^j, \Delta)$ it
is enough to prove that $\tilde H_*(\Delta,\Q) =0$ in all dimensions. $\Delta$
is the union of $\binom{j}{2}$ spaces $\Delta_{r,s} \simeq (\RP^2)^{j-1}$,
$1\leq r < s \leq j$, defined by the coincidence of the $r$-th and $s$-th
members of the sequences of $j$ points in $\RP^2$. Consider the Mayer-Vietoris
spectral sequence of this union. All (maybe, multiple) intersections of spaces
$\Delta_{r,s}$ are homeomorphic to $(\RP^2)^m$, $m \in [1, j-1]$, in particular
have the rational homology of a point. Therefore this spectral sequence reduces
to the one calculating the homology of the order complex of all such multiple
intersections. This order complex has a maximal element (corresponding to the
non-empty intersection of all these spaces), hence is contractible. \hfill
$\Box$ \medskip

\begin{corollary} \label{lemneord}
For any $j \geq 2$, the groups $\overline{H}_*(B(\RP^2,j),\Q)$ and
$\overline{H}_*(B(\RP^2,j),\pm \Q)$ are trivial in all dimensions.
\end{corollary}

Indeed, $\overline{H}_*(I(\RP^2,j),\Q) \simeq \overline{H}_*(B(\RP^2,j),
\nu_!(\Q)),$ where $\nu_!(\Q)$ is the direct image of the constant $\Q$-sheaf
on $I(\RP^2,j)$. $\nu_!(\Q)$ is a $j!$-dimensional local system on
$B(\RP^2,j)$, containing as direct summands both local systems $\Q$ and $\pm
\Q$. So, both groups mentioned in Corollary \ref{lemneord} are direct summands
of the group $\overline{H}_*(I(\RP^2,j),\Q)$, which is trivial by Lemma
\ref{lemord}. \hfill $\Box$ \medskip

Recall the notation $B^\times(\RP^2,4)$ for the subset in $B(\RP^2,4)$
consisting of collections not contained entirely in any line in $\RP^2$.

\begin{lemma} \label{lem115}
The Borel--Moore homology groups $\overline{H}_*(B^\times(\RP^2,4),\Q)$ and
$\overline{H}_*(B^\times(\RP^2,4),\pm \Q)$ are trivial in all dimensions.
\end{lemma}

{\it Proof.} Denote by $T$ the difference $B(\RP^2,4) \setminus
B^\times(\RP^2,4)$. By Corollary \ref{lemneord} and exact sequence of the pair
$(B(\RP^2,4), T)$, it is enough to prove that the groups $\overline{H}_*(T,
\Q)$ and $\overline{H}_*(T, \pm \Q)$ are trivial. By Poincar\'e duality, this
is equivalent to the triviality of groups $H_*(T, \Q \otimes {\mathbb O}r)$ and
$H_*(T, \pm \Q \otimes {\mathbb O}r)$.

$T$ is fibered over the space $\widehat{\RP^2}$ of all lines in $\RP^2$. Its
fiber over the point, corresponding to the line $l$, is the configuration space
$B(l,4)$. By Lemma \ref{lem1} this fiber is homotopy equivalent to $S^1$. The
retractions of these fibers to the circles described in Remark \ref{rem00} can
be performed uniformly, therefore $T$ contains a $3$-dimensional deformation
retract $\tilde T$ also fibered over $\RP^2$; its fiber over the projective
line $l$  consists of all quadruples of points separating this line (considered
as a circle in the unit sphere factored by the central symmetry) into arcs of
length $\pi/4$. So, it remains to prove the triviality of groups
\begin{equation}
H_*(\tilde T, \Q \otimes \widetilde{{\mathbb O}r}) \qquad \mbox{and} \qquad
H_*(\tilde T, \pm \Q \otimes \widetilde{{\mathbb O}r}), \label{both}
\end{equation} where $\widetilde {{\mathbb O}r}$ is the restriction to $\tilde
T$ of the orientation sheaf of $T$.

The group $SO(3)$ (and hence also it double cover $SU(2)$) acts transitively on
$\tilde T$; this action defines a principal $16$-fold covering $SU(2) \to
\tilde T$. The 16-dimensional direct image of the trivial $\Q$-bundle under
this covering contains both sheaves $\Q \otimes \widetilde{{\mathbb O}r})$ and
$\pm \Q \otimes \widetilde{{\mathbb O}r}$ as direct summands, hence both groups
(\ref{both}) are at most 1-dimensional in dimensions 0 and 3 and are trivial in
all other dimensions. However, both these 1-dimensional sheaves are not
equivalent to the constant one, therefore the groups (\ref{both}) are trivial
in dimension 0; by the Euler characteristic considerations the same is true for
these groups in dimension 3. \hfill $\Box$

\begin{lemma}
\label{lem11} The group $\overline{H}_i(B(\RP^2,2), {\mathbb O}r \otimes \Q)$
is equal to $\Q$ for $i=4$ and $i=1$, and is trivial in all other dimensions.
\end{lemma}

{\it Proof.} This group is Poincar\'e dual to $H_{4-i}(B(\RP^2,2), \Q)$; let us
calculate the latter group. Associating with any two-point configuration the
line spanned by these points, we obtain the fibration of $B(\RP^2,2)$ over
$\RP^2$. Its fiber is homotopy equivalent to $S^1$. The generator of
$\pi_1(\RP^2)$ acts non-trivially on the group $H_1$ of the fiber, thus the
term $E^2$ of the homological spectral sequence of this fibration has only two
non-trivial cells, $E^2_{0,0} \simeq \Q$ and $E^2_{2,1} \simeq H_2(\RP^2,
{\mathbb O}r \otimes \Q) \simeq \Q$. \hfill $\Box$

\begin{proposition}[Caratheodory theorem, see e.g. \cite{Vfil}, \cite{KK}]
\label{carat} For any $r \ge 1$, the $r$-th self-join $(S^1)^{*r}$ of $S^1$ is
homeomorphic to $S^{2r-1}$. \quad $\Box$ \end{proposition}

\begin{lemma} \label{lem33}
If $r$ is even, then the sphere $(S^1)^{*r} \sim S^{2r-1}$ has a canonical
orientation, and hence a canonical generator of its $(2r-1)$-dimensional
homology group. If $r$ is odd, then the homeomorphism of $(S^1)^{*r}$ to
itself, induced by any orientation-reversing automorphism of $S^1$, violates
the orientations of $(S^1)^{*r}$.
\end{lemma}

{\it Proof}. The group $H_{2r-1}((S^1)^{*r})$ is generated by the Borel--Moore
fundamental class of an open dense subset in $(S^1)^{*r}$, which is fibered
over the configuration space $B(S^1,r)$. Its fiber over an $r$-configuration is
the open $(r-1)$-dimensional simplex, whose vertices correspond to the points
of this configuration. The orientation of this space consists of the
orientations of the base and the fibers. Given an $r$-configuration in $S^1$,
let us number arbitrarily its points, and define the $j$-th tangent vector of
the base, $j=1, \dots, r$, as a shift of the $j$-th point along the chosen
orientation of $S^1$, keeping all remaining $r-1$ points unmoved. The fiber
over this configuration will be oriented by the same order of these $r$ points.
Renumbering these points, we simultaneously change or preserve the orientations
of the base and the fibre, thus preserving the orientation of the total space.
On the other hand, changing the chosen orientation of $S^1$ we reverse $r$
basic vectors of the orienting frame. \hfill $\Box$ \medskip

\section{Calculation of the spectral sequence}

By item 2 of Lemma \ref{lem2} and the Thom isomorphism,
we have
\begin{equation} \label{thom}
\overline{H}_m(\sigma_p \sm \sigma_{p-1}, \Q) \simeq \overline{H}_{m-d(p)\cdot
k}(\Psi_p \setminus \Psi_{p-1}, \Q \otimes \mbox{O})
\end{equation}
for any $m$ and $p$; here $\mbox{O}$ is the orientation sheaf of the vector
bundle, whose fiber over any point of $\breve \Psi(X)$, $X \in J_p$, is the
space $L(X) \subset W^k$ of all quadratic systems vanishing on the set $\{X\}$,
and $d(p)$ is the dimension of these spaces $L(X), X \in J_p$ (i.e. $d(1)=5$,
$d(2)=4$, $d(3)=d(4)=d(5)=3$, $d(6)=d(7)=2$, $d(8)=d(9)=1$, and $d(10)=0$). By
the construction, this vector bundle is the direct sum of $k$ copies of a
bundle, whose fibers are corresponding subspaces of $W$. Therefore this
orientation sheaf is trivial if $k$ is even, and coincides with the orientation
sheaf of the latter vector bundle if $k$ is odd.

So, we are going to calculate the right-hand groups in (\ref{thom}).

\subsection{First column}

The first term $\sigma_1$ of our filtration of $\sigma$ is the space of a
$5k$-dimensional vector bundle over $\RP^2$. This vector bundle is naturally
oriented for any $k$. Indeed, any of its $k$ summands is a 5-dimensional vector
subbundle of the trivial bundle $\RP^2 \times W \to \RP^2$; its fiber over a
point $x \in \RP^2$ is the space of all polynomials $f \in W^k$ vanishing on
the line $\{x\} \subset \R^3$. The canonical coorientation of such a fiber in
$W$ is defined by the increase of the polynomials $f$ on the corresponding line
$\{x\}$. By (\ref{thom}) we have
\begin{equation} \label{fsterm}
E^1_{1,q} \equiv \overline{H}_{1+q-5k}(\RP^2, \Q).
\end{equation}
So, the group $E^1_{1,q}$ is equal to $\Q$ if $q=5k-1$ and is trivial otherwise.

\subsection{Second column}
\label{2ndclmn}

\begin{lemma}
\label{lem48} The group $\overline{H}_*(\sigma_2 \setminus \sigma_1, \Q)$ is
trivial in all dimensions.
\end{lemma}

{\it Proof.} $\sigma_2 \setminus \sigma_1$ is the space of a fiber bundle over
the configuration space $B(\RP^2,2)$. Namely, it is the fibered product of

--- the $(4k)$-dimensional vector bundle, whose fiber over the point $\{x, y\}
\in B(\RP^2,2)$ consists of all systems $F \in W^k$ vanishing on both lines in
$\R^3$ corresponding to the points $x$ and $y$, and

--- the bundle of open intervals, whose endpoints are associated with the
points $x$ and $y$. Any such interval is constituted by two half-intervals
$(\{x\}, \{x \cup y\}]$ and $(\{y\}, \{x \cup y\}]$ in $\breve \Psi_2$, where
$\{x\}$ and $\{y\}$ are points of $J_1$ corresponding to the points $x$ and
$y$, and $\{x \cup y\} \in J_2$. The endpoints $\{x\}$ and $\{y\}$ are excluded
from the segment as they belong to the first term $\Psi_1$ of our filtration.

Any loop in $B(\RP^2,2)$ not permuting two points of the configuration does not
violate the orientations of these two bundles. A loop permuting these points
changes the orientations of all $k$ factors of the first bundle, as well as the
orientation of the second. Thus the $m$-dimensional Borel--Moore homology group
of the total space of this fiber bundle is equal to
$\overline{H}_{m-4k-1}(B(\RP^2, 2), \Q)$ if $k$ is odd, and to
$\overline{H}_{m-4k-1}(B(\RP^2, 2), \pm \Q)$ if $k$ is even. By Corollary
\ref{lemneord} of Lemma \ref{lemord} all these groups are trivial.  \hfill
$\Box$
\medskip

Therefore the column $E^1_{2,*}$ of the main spectral sequence consists of
zeroes only.

\subsection{Third column}
\label{3dclmn}

The space $\sigma_3 \setminus \sigma_2$ is fibered over the configuration space
$J_3 = B(\RP^2,3)$. Its fiber over a 3-point configuration $X \in J_3$ is equal
to the product of the $3k$-dimensional vector space $L(X)$ and the {\em open}
$2$-dimensional simplex (whose vertices are associated with the points of the
configuration, and the entire simplex is the order complex of non-empty subsets
of this configuration). Any element of $\pi_1( B(\RP^2,3))$, defining an odd
permutation of three points, changes the orientations of the bundle of
triangles and of all $k$ summands of the $3k$-dimensional vector bundle.
Therefore we have the following lemma.

\begin{lemma}
\label{lem35} For even $k$, $\overline{H}_m(\sigma_3 \setminus \sigma_2, \Q) =
\overline{H}_{m-3k-2}(B(\RP^2,3), \pm \Q).$ For odd $k$,
$\overline{H}_m(\sigma_3 \setminus \sigma_2, \Q) =
\overline{H}_{m-3k-2}(B(\RP^2,3), \Q).$
\end{lemma}

By Corollary \ref{lemneord} of Lemma \ref{lemord} all these groups are trivial
in all dimensions, so the column $p=3$ of the main spectral sequence is
trivial.

\subsection{Fourth column}
\label{4thclmn}

The space $\sigma_5 \setminus \sigma_3$ is fibered over the space
$\widehat{\RP^2}$ of all lines in $\RP^2$. Its fiber over the line $l$ is equal
to the product of the $3k$-dimensional subspace $L(l) \subset W^k$ (consisting
of all polynomial systems vanishing on $l$) and the space $\Psi(l) \setminus
\Psi_3$.

The space $\Psi(l)$ is a compact cone with the vertex $\{l\} \in J_5$, in
particular is contractible. Its subset $\Psi(l) \cap \Psi_3$ is the order
complex of all subsets of cardinality $\le 3$ in the line $l \sim S^1$. It is
the union of order subcomplexes $\Psi(\{a,b,c\})$, where $a, b, c$ are
arbitrary distinct points of $l$. Any such subcomplex can be identified with
the triangle with vertices $\{a\}, \{b\}$ and $\{c\}$; the entire space
$\Psi(l) \cap \Psi_3$ becomes thus identified with the union of all such
simplices supplied with the natural topology. By Proposition \ref{carat} this
space is homeomorphic to $S^5$. The exact sequence of the pair $(\Psi(l),
\Psi(l) \cap \Psi_3)$ proves now that the Borel--Moore homology of the fiber $
\Psi(l) \setminus \Psi_3$ is equal to $\Q$ in dimension 6 and is trivial in all
other dimensions. By Lemma \ref{lem33} the monodromy over $\widehat{\RP^2}$
reverses the orientation of this factor of the fiber, and also of this homology
group.

On the other hand, the bundle with fiber $L(l)$ splits into the direct sum of
$k$ 3-dimensional bundles, any of which changes its orientation over the
generator of $\pi_1(\widehat{\RP^2})$. So, we obtain the following fact.

\begin{lemma}
\label{lem39} For even $k$, $\overline{H}_i(\sigma_5 \setminus \sigma_3, \Q)
\simeq \overline{H}_{i-3k-6}(\widehat{\RP^2},\Q \otimes {\mathbb O}r)$. In
particular, it is equal to $\Q$ if $i=3k+8$ and is trivial in all other
dimensions.

For odd $k$, $\overline{H}_i(\sigma_5 \setminus \sigma_3, \Q) \simeq
\overline{H}_{i-3k-6}(\widehat{\RP^2},\Q)$. In particular, it is equal to $\Q$
if $i=3k+6$ and is trivial in all other dimensions. \hfill $\Box$
\end{lemma}

This gives us the columns $p=4$ of both parts of Fig. \ref{mss}. \hfill $\Box$

\subsection{Fifth column}

The space $\sigma_6 \setminus \sigma_5$ is fibered over the restricted
configuration space $B^\times(\RP^2,4)$. Its fiber over $X\in
B^\times(\RP^2,4)$ is the product of the $2k$-dimensional vector space $L(X)$
and the open $3$-dimensional simplex, whose vertices correspond to the points
of the collection $\{X\}$. In exactly the same way as in \S\S \ref{2ndclmn},
\ref{3dclmn}, we obtain the following fact.

\begin{lemma} \label{lem4a1} The group
$\overline{H}_m(\sigma_6 \setminus \sigma_5, \Q)$ is equal to
$\overline{H}_{m-2k-3}(B^\times(\RP^2,4), \Q)$ if $k$ is odd, and to
$\overline{H}_{m-2k-3}(B^\times(\RP^2,4), \pm \Q)$ if $k$ is even.
\end{lemma}

By Lemma \ref{lem115} all these groups are trivial for any $m$, and the column
$p=5$ of our spectral sequence is zero.

\subsection{Sixth column}

\begin{lemma} \label{lem4b1}
For odd $k$, $\overline{H}_i(\sigma_7 \setminus \sigma_6,\Q) \simeq
\overline{H}_{i-2k-9}(\RP^2, \Q \otimes {\mathbb O}r)$, in particular this
group is equal to $\Q$ if $i=2k+11$ and is trivial otherwise. For even $k$,
$\overline{H}_i(\sigma_7 \setminus \sigma_6,\Q) \simeq
\overline{H}_{i-2k-9}(\RP^2, \Q)$, in particular is equal to $\Q$ if $i=2k+9$
and is trivial otherwise.
\end{lemma}

{\it Proof.} $\sigma_7 \setminus \sigma_6$ is the space of a fiber bundle over
the space $J_7$ of configurations
\begin{equation}
\label{conf4b} \{\mbox{a line } \{l\} \in \widehat{\RP^2}; \mbox{a point } u
\in \RP^2 \setminus l\}.\end{equation} Its fiber over such a configuration
$X=\{l,u\} \in J_7$ is the product
\begin{equation} \label{fiber42}
\breve \Psi(X) \times L(X),
\end{equation}
where $L(X) \subset W^k$ is the $2k$-dimensional space of all systems $F \in
W^k$ vanishing on $l$ and on $u$. The space of configurations (\ref{conf4b}) is
obviously diffeomorphic to $T_* \RP^2$, in particular is orientable.

\begin{lemma}
\label{lem417} 1. The non-trivial element of the fundamental group
$\pi_1(T_*\RP^2) \simeq \pi_1(\RP^2)$ of this configuration space violates the
orientation of any of $k$ two-dimensional summands of the second factor $L(X)$
of the fiber $($\ref{fiber42}$)$.

2. The link $\partial \Psi(X)$ of the first factor of the fiber
$($\ref{fiber42}$)$ is homotopy equivalent to $S^6$, in particular the group
$\overline{H}_i(\breve \Psi(X), \Q)$ is isomorphic to $\Q$ if $i=7$ and is
trivial for all other $i$.

3. The monodromy over the non-trivial element of $\pi_1(T_*\RP^2)$ acts as
multiplication by $-1$ on the homology bundle associated with the latter
factors.
\end{lemma}

{\it Proof.} 1. For the loop in our configuration space, generating its
fundamental group, we choose the family of pairs $\{l_\varphi, u_\varphi\}$,
$\varphi \in [0, \pi]$, where the line $l_\varphi$ is defined by the equation
$x \cos \varphi + y \sin \varphi =0$, and the point $u_\varphi$ by homogeneous
coordinates $(\cos \varphi : \sin \varphi : 0)$. For any configuration $X =
\{l,u\}$, the fiber of any of our $k$ summands over $X$ consists of all
quadratic forms vanishing on $l$ and $u$. For two independent sections of the
two-dimensional bundle of such forms over the set of pairs $\{l_\varphi,
u_\varphi\}$, where $\varphi$ runs over the half-interval $ \varphi \in
[0,\pi)$, we choose the families $$ (x \cos \varphi + y \sin \varphi)(x \sin
\varphi - y \cos \varphi) \mbox{ \quad and \quad } (x \cos \varphi + y \sin
\varphi)z .$$ The first of these sections extends to a continuous section over
the entire closed loop, while the limit value of the second one at $\varphi =
\pi$ is opposite to the value at $\varphi =0$. So, the monodromy over this loop
acts
on the fiber as the operator {\footnotesize $\left| \begin{array}{cc} 1 & 0 \\
0 & -1 \end{array} \right|$}.

2. The link $\partial \Psi(X),$ $X=\{l,u\}$, consists of two subspaces: $A=$
the order complex $\Psi(l),$ and $B=$ the union of all order complexes
$\Psi(\{a,b,c,u\})$, where $a, b$ and $c$ are arbitrary distinct points in the
line $l$. The intersection $A \cap B$ is the space $\Psi(l) \cap \Psi_3$, i.e.
the union of similar order subcomplexes $\Psi(\{a,b,c\})$. Since \S
\ref{4thclmn} we know that this intersection $A \cap B$ is homeomorphic to
$S^5$. The set $A$ is compact and contractible, hence the pair $(A, A \cap B)$
is homotopy equivalent to the pair $(\mbox{cone over } A \cap B, \mbox{its base
} A \cap B)$. As for the space $B$, any of order complexes $\Psi(\{a, b, c,
u\})$ constituting it is canonically PL-homeomorphic to the tetrahedron with
vertices $a, b, c$ and $u,$ and hence to the cone over $\Psi(\{a, b, c\})$ with
vertex $\{u\}$. These homeomorphisms over all triples $\{a, b, c\}$ define a
homeomorphism between the union $B$ of all these tetrahedra and another cone
over $\Psi(l) \cap \Psi_3$. Therefore our space $A \cup B$ is homotopy
equivalent to the suspension over $A \cap B \sim S^5$.

3. By Lemma \ref{lem33}, the generator of $\pi_1(\RP^2)$ acts on the group
$H_5(\Psi(l) \cap \Psi_3, \Q)$ as multiplication by $-1$. Two parts $A$ and $B$
of $\partial \Psi(\{l,u\}) \sim \Sigma (\Psi(l) \cap \Psi_3)$ are of different
origin and cannot be permuted by the monodromy, hence the group $H_6(\partial
\Psi(\{l,u\}), \Q)$ also is reversed by this generator. \hfill $\Box$
\medskip

Lemma \ref{lem417} is completely proved, and Lemma \ref{lem4b1} follows
immediately from it by a) the spectral sequence of the fiber bundle
(\ref{fiber42}) over the space $J_6 \sim T_*\RP^2$ of configurations
(\ref{conf4b}), and b) the Thom isomorphism for the non-orientable fibration of
this space over $\RP^2$. \hfill $\Box$ \medskip

This gives us columns $p=6$ of both spectral sequences of Fig. \ref{mss}.

\subsection{Seventh column}

\begin{proposition}
\label{prop5b} For any even $k$, the group $\overline{H}_*(\sigma_8 \setminus
\sigma_7,\Q)$ is trivial in all dimensions. If $k$ is odd, then
$\overline{H}_i(\sigma_8 \setminus \sigma_7,\Q)$ is isomorphic to $\Q$ for
$i=k+12$ and $i=k+9$, and is trivial for all other values of $i$.
\end{proposition}

{\it Proof.} The space $\sigma_8 \setminus \sigma_7$ is fibered over the space
$J_8 = B(\widehat{\RP^2},2)$ of 2-point configurations in the dual projective
space; its fiber over any point $X$ of $J_8$ is homeomorphic to $\breve \Psi(X)
\times \R^k$.

\begin{lemma} \label{lem5b}
For any point $X \in B(\widehat{\RP^2},2)$, the link $\partial \Psi(X)$ of the
order complex $\Psi(X)$ is homology equivalent to the sphere $S^7$.
\end{lemma}

{\it Proof.} Denote by $l$ and $l'$ two lines constituting the configuration
$X$, and by $O$ their common point. Let us filter the link $\partial \Psi(X)$:
\begin{equation} \label{fltr7}(\Psi(l) \cup \Psi(l')) \subset
(\partial \Psi(X) \cap \Psi_6) \subset \partial \Psi(X),
\end{equation}
and study this filtration.

The order complexes $\Psi(l)$ and $\Psi(l')$ are contractible by definition,
and their intersection set is a single point in $\Psi_1$, corresponding to the
point $l \cap l'$. Therefore the union $\Psi(l) \cup \Psi(l')$ is contractible
too.

The space $(\Psi(X) \cap \Psi_6) \setminus (\Psi(l) \cup \Psi(l'))$ consists of
two similar disjoint pieces, swept out by all order complexes $\Psi(\{l,a'\})
\setminus \Psi(l)$, (respectively,  $\Psi(\{l',a\}) \setminus \Psi(l')$, )
where $a'$ runs over the affine line $l' \setminus O$ and $a$ runs over $l
\setminus O$. Obviously, the first piece is homeomorphic to the direct product
$(l' \setminus O) \times (\Psi(\{l,a'\}) \setminus \Psi(l))$ for an arbitrary
$a' \in l' \setminus O$. By the proof of Lemma \ref{lem417}(2), the pair
$(\Psi(\{l,a'\}) , \Psi(l))$ is homotopy equivalent to the pair $(D^7, D^6)$,
where $D^7$ is the 7-dimensional ball, and $D^6$ is a closed hemisphere of its
boundary $\partial D^7 \sim S^6$. In particular, the group $H_m(\Psi(\{l,a'\}),
\Psi(l)) \equiv \overline{H}_m(\Psi(\{l,a'\}) \setminus \Psi(l))$ is trivial
for any $m$. Hence by the K\"unneth formula all Borel-Moore homology groups of
the first piece of $(\Psi(X) \cap \Psi_6) \setminus (\Psi(l) \cup \Psi(l'))$
also are trivial. The second piece is homeomorphic to it. So, the entire second
term of the filtration (\ref{fltr7}) also has the homology of a point.

Finally, the space $\partial \Psi(X) \setminus \Psi_6$ is the union of all
spaces $\breve \Psi(\{a,b,c,d\})$ where $a \neq b \in l \setminus O$, and $c
\neq d \in l' \setminus O$. These spaces are homeomorphic to an open
tetrahedron, and the base of this family, $B(\R^1,2) \times B(\R^1, 2)$ is
homeomorphic to the 4-dimensional open ball. In particular, the Borel--Moore
homology group of this space is equal to $\Q$ in dimension $7$ and trivial in
all other dimensions. \hfill $\Box$
\medskip

Now, let us study the action of $\pi_1(B(\widehat{\RP^2},2))$ on the
Borel--Moore homology group of fibers $\breve \Psi(X) \times L(X)$ of our
fibration of $\sigma_7 \setminus \sigma_6$ over $B(\widehat{\RP^2},2)$. By the
previous lemma this group is equal to $\Q$ in dimension $8+k$ and is trivial in
all other dimensions.

Let us choose for the basepoint in $B(\widehat{\RP^2},2)$ the pair of lines,
one of which is obtained from meridians with longitudes 0 and 180, and the
other with longitudes 90 and 270. The group $\pi_1(B(\widehat{\RP^2},2))$ is
generated by two elements, $\theta$ and $\varpi$. Namely, $\theta$ is defined
by the $90^\circ$ rotation of the sphere around the main axis containing the
poles; $\varpi$ rotates the first line by $180^\circ$ around the axis
connecting the equator point with coordinates $(0,0)$ and $(0,180)$, and keeps
the second line unmoved.

\begin{lemma} \label{lem5b2}
Both loops $\theta, \varpi \in \pi_1(B(\widehat{\RP^2},2))$ change the
orientation of $B(\widehat{\RP^2},2)$ and of any of $k$ factors of the fiber
bundle with fibers $L(X)$. The one-dimensional homology groups
$\overline{H}_8(\breve \Psi(X),\Q)$ and $H_7(\partial \Psi(X),\Q)$ have a
canonical orientation for all $X$, depending continuously on $X$; in particular
the group $\pi_1(B(\widehat{\RP^2},2))$ acts trivially on them.
\end{lemma}

{\it Proof.} The statement on the orientation of the base follows immediately
from Lemma \ref{lem77}. Any of $k$ summands of the fiber $L(X)$ is the space of
quadratic functions vanishing on our pair of lines. We can choose one of two
domains into which these lines divide $\RP^2$, and orient any summand by the
increase of our quadratic functions inside this domain. The continuation along
any of loops $\theta$ or $\varpi$ permutes these two domains, hence violates
this orientation.

Now we define an invariant orientation of the homology bundle whose fiber over
$X \in B(\widehat{\RP^2},2)$ is the group $H_7(\partial \Psi(X),\Q)$. This
group is generated by the fundamental class of the 7-dimensional cell described
in the last paragraph of the proof of Lemma \ref{lem5b}. To orient this cell,
we need to choose the orientations of its base $B(l \setminus O, 2) \times B(l'
\setminus O, 2)$ and, for any point of this base (i.e. a configuration of four
points), of the fiber over it, i.e. the 3-dimensional simplex, whose four
vertices correspond to the points of this configuration. To do it, we fix an
arbitrary numbering of these four points, and define the $j$th basic tangent
vector to the base, $j =1, \dots, 4$, as a shift of the $j$th point in the
direction from the other point of the configuration in the same line and
towards the intersection point $O$ of these lines. The orientation of the fiber
(i.e. a tetrahedron) is defined by the same numbering of its vertices. A
different numbering of four points simultaneously changes or preserves the
orientations of the base and the fiber, and hence preserves the orientation of
the total space.

The orientation thus defined is obviously invariant and defines a
trivialization of the homology bundle over $B(\widehat{\RP^2},2)$ with the
fiber $H_7(\partial \Psi(X),\Q)$ or, which is the same by the boundary
isomorphism, of the group $\overline{H}_8(\breve \Psi(X),\Q)$. \hfill $\Box$
\medskip

\begin{corollary} \label{cor5b}
If $k$ is even $($respectively, odd$)$, then $\overline{H}_i(\sigma_7 \setminus
\sigma_6,\Q)$ is equal to $\overline{H}_{i-k-8}(B(\widehat{\RP^2},2), \Q)$
$($respectively, $\overline{H}_{i-k-8}(B(\widehat{\RP^2},2), {\mathbb O}r
\otimes \Q))$ for any $i$. \end{corollary}

This fact follows from Lemma \ref{lem5b2} by the Thom isomorphism. \hfill
$\Box$
\medskip

Proposition \ref{prop5b} follows from it by Corollary \ref{lemneord} and Lemma
\ref{lem11}, and justifies the columns $p=7$ in both parts of Fig. \ref{mss}.
\hfill $\Box$

\subsection{Eighth column}

\begin{lemma} \label{lemj8}
The space $J_9$ of non-empty non-singular conics in $\R^3$ is diffeomorphic to
the space of a 3-dimensional vector bundle over $\RP^2$, namely of the
symmetric square of the cotangent bundle $S^2 T^*(\RP^2)$. In particular, it is
an orientable 5-dimensional manifold homotopy equivalent to $\RP^2$.
\end{lemma}

{\it Proof.} Any such conic can be specified uniquely (up to the positive
scalings) by a quadratic form of signature $(1,2)$. Let us fix an arbitrary
Euclidean structure in $\R^3$. Then the line of positive eigenvectors of these
quadratic forms indicates a point of $\RP^2$, and the plane of negative ones
can be naturally identified with the tangent plane of $\RP^2$ at this point.
Let us always choose our quadratic form in such a way that it takes value 1 on
its positive eigenvectors of length 1. Then we obtain an isomorphism between
the space of all conics of this signature with this positive eigenvector on one
hand, and the space of negative quadratic forms on the tangent space of $\RP^2$
at this point $u$ on the other. The latter space is a convex open subset in the
space $S^2 T_u^*(\RP^2)$, hence is diffeomorphic to the entire this space (and
these diffeomorphisms can be performed uniformly over all fibers). \hfill
$\Box$

\begin{lemma} \label{lem5a}
For arbitrary $k$, the group $\overline{H}_i(\sigma_9 \setminus \sigma_8, \Q)$
is isomorphic to $\overline{H}_{i-k-11}(\RP^2, \Q \otimes {\mathbb O}r)$, i.e.
to $\Q$ if $i=k+13$ and to the trivial group in all other dimensions.
\end{lemma}

{\it Proof.} The space $\sigma_9 \setminus \sigma_8$ is fibered over $J_9$; its
fiber over a conic $C \in J_9$ is equal to $\breve \Psi(C) \times L(C)$, $L(C)
\simeq \R^k$.

The spaces $L(C)$ form a trivial vector bundle over $J_9$. Indeed, any of its
one-dimensional summands can be canonically oriented, because the space of
quadratic forms vanishing on a fixed non-empty cone is oriented in an invariant
way (say, from quadratic forms of signature $(2,1)$ to those of signature
$(1,2)$).

By Proposition \ref{carat}, the link $\partial \Psi(C)$ is homeomorphic to
$S^7$. In particular, the Borel--Moore homology groups of the fibers $\breve
\Psi(C)$ are isomorphic to $\Q$ in dimension $8$, and are trivial in all other
dimensions. By Lemma \ref{lem33} the fundamental group of the base $J_9$ acts
trivially on these homology groups.

Now, our lemma follows from the Thom isomorphism $\overline{H}_i(\sigma_9
\setminus \sigma_8,\Q) \simeq \overline{H}_{i-k-8}(J_9, \Q)$ of our orientable
fiber bundle over the space of conics, and the Thom isomorphism
$\overline{H}_j(J_9, \Q) \simeq H_{j-3}(\RP^2, \Q \otimes {\mathbb O}r)$ for
the fibration of $J_9$ over $\RP^2$, see Lemma \ref{lemj8}. \hfill $\Box$
\medskip

This proves the columns $p=8$ of both tables in Fig. \ref{mss}. In particular,
we see that in the case of odd $k$ both groups $E^1_{8, 5} \equiv
\overline{H}_{k+13}(\sigma_9 \setminus \sigma_8,\Q)$ and $E^1_{7,5} \equiv
\overline{H}_{k+12}(\sigma_8 \setminus \sigma_7,\Q)$ are equal to $\Q$.

\subsection{The differential $d_1: E^1_{8,k+5} \to E^1_{7,k+5}$ for odd $k$}

\begin{proposition} \label{lem55s}
The boundary homomorphism $$\overline{H}_{k+13}(\sigma_9 \setminus \sigma_8,\Q)
\stackrel{\partial}{\to} \overline{H}_{k+12}(\sigma_8 \setminus \sigma_7,\Q)$$
from the exact sequence of the triple is an isomorphism for any odd $k$.
\end{proposition}

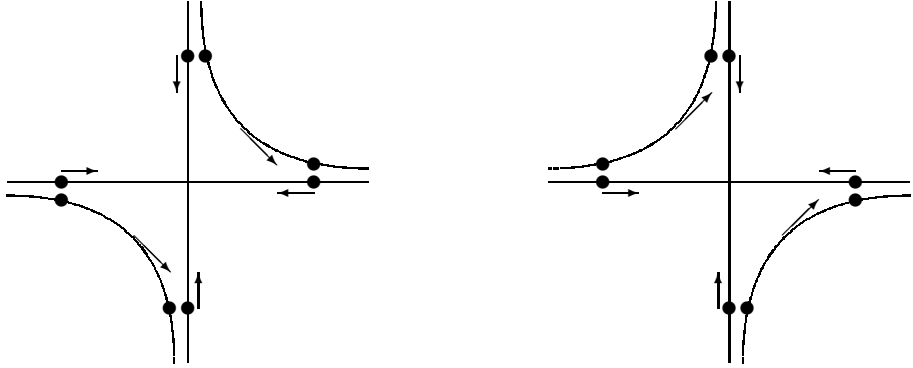
\begin{figure}
\unitlength=1.2mm
\begin{picture}(150,40)
\put(0,20){\line(1,0){40}} \put(20,0){\line(0,1){40}}
\put(60,20){\line(1,0){40}} \put(80,0){\line(0,1){40}}
\bezier{140}(0,18.5)(18.5,18.5)(18.5,0)
\bezier{140}(21.5,40)(21.5,21.5)(40,21.5)
\bezier{140}(60,21.5)(78.5,21.5)(78.5,40)
\bezier{140}(81.5,0)(81.5,18.5)(100,18.5)

\put(6,20){\circle*{1.5}} \put(34,20){\circle*{1.5}} \put(66,20){\circle*{1.5}}
\put(94,20){\circle*{1.5}} \put(20,6){\circle*{1.5}} \put(20,34){\circle*{1.5}}
\put(80,6){\circle*{1.5}} \put(80,34){\circle*{1.5}} \put(6,18){\circle*{1.5}}
\put(34,22){\circle*{1.5}} \put(66,22){\circle*{1.5}}
\put(94,18){\circle*{1.5}} \put(18,6){\circle*{1.5}} \put(22,34){\circle*{1.5}}
\put(82,6){\circle*{1.5}} \put(78,34){\circle*{1.5}}
\put(6,21.2){\vector(1,0){4}} \put(34,18.8){\vector(-1,0){4}}
\put(18.8,34){\vector(0,-1){4}} \put(21.2,6){\vector(0,1){4}}
\put(66,18.8){\vector(1,0){4}} \put(94,21.2){\vector(-1,0){4}}
\put(81.2,34){\vector(0,-1){4}} \put(78.8,6){\vector(0,1){4}}
\put(14.1,14.1){\vector(1,-1){4}} \put(25.9,25.9){\vector(1,-1){4}}
\put(74.1,25.9){\vector(1,1){4}} \put(85.9,14.1){\vector(1,1){4}}
\end{picture}
\caption{Different approximations of $J_7$ from $J_8$} \label{5atob}
\end{figure}

{\it Proof.} Let $X$ be an arbitrary point of the stratum $J_8$, i.e. a pair of
crossing lines in $\RP^2$. This stratum has codimension 1 in the space of
conics, and the open stratum $J_9$ approaches it at $X$ from two different
sides, see Fig. \ref{5atob}. Recall that $\sigma_9 \setminus \sigma_8$ is
fibered over $J_9$ with the fiber over $C \in J_9$ equal to $\breve \Psi(C)
\times L(C)$, $L(C) \simeq \R^k$. The base $J_9$ is oriented as an open
subvariety of the space $\RP^5$ of all conics, the canonical orientation of the
fibers $L(C)$ (and moreover of any of their $k$ summands) was specified in the
proof of Lemma \ref{lem5a}, and the canonical orientation of fibers $\breve
\Psi(C)$ was defined in the proof of Lemma \ref{lem33}. In a neighborhood $U$
of the point $X$, the manifold $J_8$ can be considered as the boundary of any
of two local components, into which it separates $J_9$; therefore the
restriction of the orientation of $J_9$ to any of these components induces the
boundary orientation on $J_8 \cap U$.

\begin{lemma} \label{lem5or}
When the points $C \in J_9$ tend from two different local components of $J_9$
to a point $X$ of the hypersurface $J_8 \subset \RP^5$,

1. Two boundary orientations of $J_8$ at $X$, induced from the restrictions of
one and the same orientation of $J_9$ to these components, are opposite to one
another;

2. The 7-dimensional fundamental classes of canonically oriented fibers
$\partial \Psi(C) \sim S^7$ tend to the fundamental class of the fiber
$\partial \Psi(X) \setminus \Psi_6$ supplied with the canonical orientation
$($specified in the last paragraph of the proof of Lemma \ref{lem5b}$)$ and
generating the Borel--Moore homology group of entire $\partial \Psi(X)$.

3. The canonical orientations of any of $k$ summands of the $k$-dimensional
vector bundle with fibers $L(C)$ tend to opposite orientations of the similar
line over $X$.
\end{lemma}

By the boundary isomorphisms, item 2 of Lemma \ref{lem5or} is equivalent to the
following statement: the 8-dimensional fundamental classes of canonically
oriented fibers $\breve \Psi(C)$ (homeomorphic to open 8-dimensional balls)
tend to the canonically oriented fundamental class of the piece of $\breve
\Psi(X)$ generating its 8-dimensional homology group.
\medskip

{\it Proof of Lemma \ref{lem5or}}. 1. Statement 1 follows from the fact that
the union of $J_9$ and $J_8$ is an open subset in the orientable manifold
$\RP^5$ of all conics.

2. These orientations, as they were defined, depend only on the directions of
the shifts of the points of 4-configurations. For $J_9$, it are the directions
following one and the same orientation of our conic, for $J_8$ they go towards
the crossing point, see Fig. \ref{5atob}. The limit positions of the former
directions, when conics degenerate to a cross, differ from the latter ones an
even number of times, therefore the corresponding orientations do coincide.

3. For any point $C \in J_9$ any such 1-dimensional summand consists of all
quadratic functions vanishing on the corresponding cone. The basic vector,
defining the canonical orientation of this line, is a quadratic function
positive inside the simply-connected component of the complement of the cone in
$\R^3$. The limits of such functions for the left- and right-hand parts of Fig.
\ref{5atob} are functions vanishing on the coordinate cross, but having
opposite signs outside it; therefore they define opposite orientations of the
corresponding summand of the fiber $L(X)$ over the point $X \in J_8$. \hfill
$\Box$ \medskip

So, in total the boundary orientations of $\sigma_8 \setminus \sigma_7$,
induced from these two sides, differ by the factor $(-1) \times 1 \times
(-1)^k$, which in the case of odd $k$ is equal to 1 and proves the
non-triviality of the boundary operator, i.e. the assertion of Proposition
\ref{lem55s}. \hfill $\Box$

\subsection{Ninth and last column}

\begin{proposition} \label{prop6m}
The group $\overline{H}_i(\sigma_{10} \setminus \sigma_9, \Q)$ is equal to $\Q$
if $i=14$ and is trivial in all other dimensions.
\end{proposition}

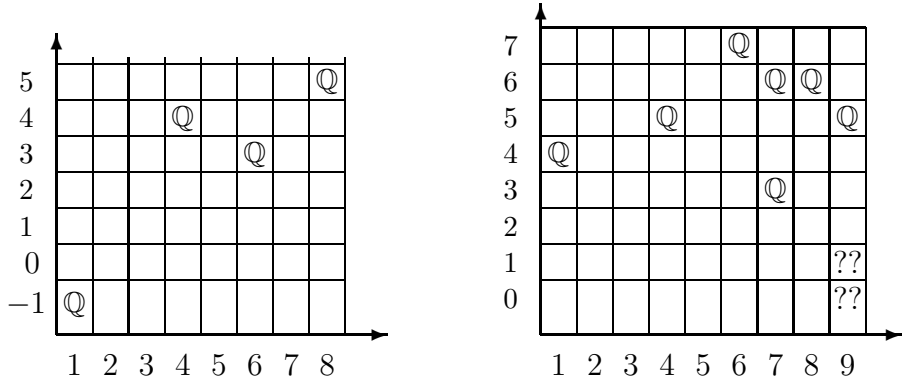
\begin{figure}
\unitlength=0.80mm \special{em:linewidth 0.4pt} \linethickness{0.4pt} \mbox{
\begin{picture}(65.00,60.00)
\thicklines \put(10.00,10.00){\vector(1,0){55.00}}
\put(10.00,10.00){\vector(0,1){50.00}} \thinlines
\put(16.00,10.00){\line(0,1){46.00}} \put(22.00,10.00){\line(0,1){46.00}}
\put(28.00,10.00){\line(0,1){46.00}} \put(34.00,10.00){\line(0,1){46.00}}
\put(40.00,10.00){\line(0,1){46.00}} \put(46.00,10.00){\line(0,1){46.00}}
\put(52.00,10.00){\line(0,1){46.00}} \put(58.00,10.00){\line(0,1){46.00}}
\put(10.00,19.00){\line(1,0){48.00}} \put(10.00,25.00){\line(1,0){48.00}}
\put(10.00,31.00){\line(1,0){48.00}} \put(10.00,37.00){\line(1,0){48.00}}
\put(10.00,43.00){\line(1,0){48.00}} \put(10.00,49.00){\line(1,0){48.00}}
\put(10.00,55.00){\line(1,0){48.00}} \put(5.00,15.00){\makebox(0,0)[cc]{$-1$}}
\put(6.00,22.00){\makebox(0,0)[cc]{$0$}}
\put(13.00,15.00){\makebox(0,0)[cc]{\small $\Q$}}
\put(5.00,46.00){\makebox(0,0)[cc]{\small $4$}}
\put(5.00,40.00){\makebox(0,0)[cc]{\small $3$}}
\put(5.00,34.00){\makebox(0,0)[cc]{\small $2$}}
\put(5.00,28.00){\makebox(0,0)[cc]{\small $1$}}
\put(31.00,46.00){\makebox(0,0)[cc]{\small $\Q$}}
\put(43.00,40.00){\makebox(0,0)[cc]{\small $\Q$}}
\put(5.00,52.00){\makebox(0,0)[cc]{\small $5$}}
\put(55.00,52.00){\makebox(0,0)[cc]{\small $\Q$}}
\put(13.00,5.00){\makebox(0,0)[cc]{1}} \put(19.00,5.00){\makebox(0,0)[cc]{2}}
\put(25.00,5.00){\makebox(0,0)[cc]{3}} \put(31.00,5.00){\makebox(0,0)[cc]{4}}
\put(37.00,5.00){\makebox(0,0)[cc]{5}} \put(43.00,5.00){\makebox(0,0)[cc]{6}}
\put(49.00,5.00){\makebox(0,0)[cc]{7}} \put(55.00,5.00){\makebox(0,0)[cc]{8}}
\end{picture}
} \qquad \mbox{
\begin{picture}(70.00,65.00)
\thicklines \put(10.00,10.00){\vector(1,0){60.00}}
\put(10.00,10.00){\vector(0,1){55.00}} \thinlines
\put(16.00,10.00){\line(0,1){51.00}} \put(22.00,10.00){\line(0,1){51.00}}
\put(28.00,10.00){\line(0,1){51.00}} \put(34.00,10.00){\line(0,1){51.00}}
\put(40.00,10.00){\line(0,1){51.00}} \put(46.00,10.00){\line(0,1){51.00}}
\put(52.00,10.00){\line(0,1){51.00}} \put(58.00,10.00){\line(0,1){51.00}}
\put(64.00,10.00){\line(0,1){51.00}}

\put(10.00,19.00){\line(1,0){54.00}} \put(10.00,25.00){\line(1,0){54.00}}
\put(10.00,31.00){\line(1,0){54.00}} \put(10.00,37.00){\line(1,0){54.00}}
\put(10.00,43.00){\line(1,0){54.00}} \put(10.00,49.00){\line(1,0){54.00}}
\put(10.00,55.00){\line(1,0){54.00}} \put(10.00,61.00){\line(1,0){54.00}}

\put(13.00,5.00){\makebox(0,0)[cc]{1}} \put(19.00,5.00){\makebox(0,0)[cc]{2}}
\put(25.00,5.00){\makebox(0,0)[cc]{3}} \put(31.00,5.00){\makebox(0,0)[cc]{4}}
\put(37.00,5.00){\makebox(0,0)[cc]{5}} \put(43.00,5.00){\makebox(0,0)[cc]{6}}
\put(49.00,5.00){\makebox(0,0)[cc]{7}} \put(55.00,5.00){\makebox(0,0)[cc]{8}}
\put(61.00,5.00){\makebox(0,0)[cc]{9}}
\put(13.00,40.00){\makebox(0,0)[cc]{\small $\Q$}}
\put(31.00,46.00){\makebox(0,0)[cc]{\small $\Q$}}
\put(43.00,58.00){\makebox(0,0)[cc]{\small $\Q$}}
\put(55.00,52.00){\makebox(0,0)[cc]{\small $\Q$}}
\put(49.00,52.00){\makebox(0,0)[cc]{\small $\Q$}}
\put(61.00,46.00){\makebox(0,0)[cc]{\small $\Q$}}
\put(49.00,34.00){\makebox(0,0)[cc]{\small $\Q$}}
\put(5.00,16.00){\makebox(0,0)[cc]{\small $0$}}
\put(5.00,22.00){\makebox(0,0)[cc]{\small $1$}}
\put(5.00,28.00){\makebox(0,0)[cc]{\small $2$}}
\put(5.00,34.00){\makebox(0,0)[cc]{\small $3$}}
\put(5.00,40.00){\makebox(0,0)[cc]{\small $4$}}
\put(5.00,46.00){\makebox(0,0)[cc]{\small $5$}}
\put(5.00,52.00){\makebox(0,0)[cc]{\small $6$}}
\put(5.00,58.00){\makebox(0,0)[cc]{\small $7$}}
\put(61.00,22.00){\makebox(0,0)[cc]{$??$}}
\put(61.00,16.00){\makebox(0,0)[cc]{$??$}}
\end{picture}
} \caption{$E^1$ for $k=0$ \hspace{2cm} and $k=1$ \hspace*{2cm} }
\label{zeross}
\end{figure}

{\it Proof.} The space $J_{10}$ consists of unique point $[\RP^2]$ specified by
the zero quadratic form. Therefore $\sigma_{10} \setminus \sigma_9 = \breve
\Psi([\RP^2])$. The order complex $\Psi([\RP^2])$ can be calculated as the
space $\sigma$ of our simplicial resolution for the case $k=0$ (i.e., with
0-dimensional spaces $L(X)$). The Borel--Moore homology group of its link
$\partial \Psi([\RP^2]) \subset \sigma_9$ can thus be calculated by the
spectral sequence, whose term $E_1$ coincides with columns $ p \in [1,8]$ of
our main spectral sequence for $k=0$. These columns are calculated in the
previous sections, and the spectral sequence looks as is shown in the left-hand
part of Fig. \ref{zeross}.

\begin{lemma}
\label{lem61} The differential $\partial_2: E^2_{6,3} \to E^2_{4,4}$ of the
spectral sequence of Fig. \ref{zeross} $($left$)$ is an isomorphism.
\end{lemma}

{\it Proof.} Consider our main spectral sequence for $k=1$. All its columns
with $p=1, \dots, 8$ are already calculated, therefore it looks as is shown in
the right-hand part of Fig. \ref{zeross}, where both signs ?? in the ninth
column should be replaced by $\Q$ or $0$ depending on whether the differential
in question is trivial or not. However, the answer in the case $k=1$ is known:
the space $W^1 \setminus \Sigma$ consists of two contractible components.
Therefore by Alexander duality the group $\overline{H}_i(\Sigma, \Q) \equiv
\overline{H}_i(\sigma, \Q)$ is equal to $\Q$ for $i=5$ and is trivial in all
other dimensions. Thus all cells of our spectral sequence except for $E_{1,4}$
should be killed by some its differential. This is impossible if the group
$E^1_{9,0}$ is non-trivial. \hfill $\Box$
\medskip

So, $H_i\left(\partial \Psi([\RP^2]),\Q \right)$ is equal to $\Q$ for $i=0$ and
$i=13$, and is trivial in all other dimensions. Proposition \ref{prop6m}
follows from this fact by the exact sequence of the pair $\left(\Psi([\RP^2]),
\partial \Psi([\RP^2])\right)$. \hfill $\Box$

Theorem \ref{prop21} is completely proved. \hfill $\Box$

\section{Concluding remarks}

{\bf 1. A problem:} is there a more direct proof of Theorem \ref{mthm12}, in
the same way as Proposition \ref{toy12} can be proved by the methods of
\cite{Borel} type?

Of course, the space $W^k \setminus \Sigma$ is fibered over $S^{k-1}$: with any
non-resultant system we can associate the image of a fixed vector, say $(1,
0,0)$. However, it seems that this structure helps us quite little. The
Poincar\'e polynomials of the fibers of these fiber bundles follow easily from
Theorem \ref{mthm12}: in the case of odd (respectively, even) $k \geq 3$ they
are equal to $1 + t^{3k-7}(1+t^{k-5})(1+t^{k-2})$ (respectively, $1+t^{3k-9} +
t^{5k-14}$).
\medskip

{\bf 2.} The construction of resolutions of discriminants in the terms of
self-joins was formulated in full generality (sufficient for our current needs)
in \cite{Gorinov} (although its particular versions appeared in less difficult
situations also in \cite{VO}, \cite{V}, \cite{hompol}). It may seem that this
construction is less canonical than the one in terms of Hilbert schemes (see
e.g. \cite{novikov}). However in practice, when one needs not only to prove
abstract theorems, but to accomplish a difficult explicit homological
calculation, the approach of \cite{Gorinov} proves to be more effective.
\medskip

{\bf 3. The formal analog.} The space $W^k$ of quadratic maps $\R^3 \to \R^k$
is a subspace of the (infinite-dimensional) space of all maps $\R^3 \to \R^k$,
quadratic in restriction to any one-dimensional subspace in $\R^3$. The latter
space of maps also contains a resultant variety: it consists of maps having
common zeros in $\R^3 \setminus 0$. The corresponding space of non-resultant
maps is homotopy equivalent to the space of continuous maps $\RP^2 \to
S^{k-1}$. Its rational homology group for $k > 3$ can be easily calculated
(e.g. by the methods of \cite{V}) and coincides with that of $S^{k-1}$: all
other rational homology groups are trivial by Corollary \ref{lemneord} of the
present article. So it is much poorer than its nonstable analog calculated
above.

\end{document}